\documentclass{svjour3}
\usepackage{amssymb}
\usepackage[hidelinks]{hyperref}

\usepackage{epic}
\usepackage{amsmath}
\usepackage{amsfonts}
\usepackage{pstricks}

\usepackage{float}
\usepackage{subfig}
\usepackage{graphicx}
\usepackage{multirow}
\usepackage{adjustbox}
\usepackage{color}
\usepackage[subfigure]{tocloft}
\usepackage{subfig}
\usepackage[utf8]{inputenc}
\usepackage{soul}

\allowdisplaybreaks
\numberwithin{equation}{section}
\def\bigO{\mathcal{O}}

\begin{document}

\title{On approximate implicit Taylor methods for ordinary differential equations} 

\author{Antonio~Baeza \and Raimund~B\"{u}rger \and
  Maria del Carmen~Mart\'{\i} \and Pep~Mulet \and David~Zor\'{\i}o}
\authorrunning{A.\ Baeza, R.\ B\"{u}rger,  M.C. Mart\'{\i}, P.\ Mulet and D.\ Zor\'{\i}o}
\institute{A.\ Baeza \and M.C. Mart\'{\i} \and P.\ Mulet  \at 
  Departament de Matem\`{a}tiques \\  Universitat de
  Val\`{e}ncia \\ E-46100 Burjassot,      Spain \\ 
  \email{antonio.baeza@uv.es, Maria.C.Marti@uv.es, pep.mulet@uv.es} 
  \and
  R.\ B\"{u}rger \at  
  CI$^2$MA \& Departamento de Ingenier\'{\i}a Matem\'{a}tica \\ 
   Universidad de Concepci\'{o}n  \\ 
   Casilla 160-C,  Concepci\'{o}n,    Chile\\ 
   \email{rburger@ing-mat.udec.cl} 
   \and 
   D.\ Zor\'{\i}o \at 
    CI$^2$MA,  
   Universidad de Concepci\'{o}n  \\ 
   Casilla 160-C,  Concepci\'{o}n,    Chile\\ 
   \email{dzorio@ci2ma.udec.cl}  
}

\thispagestyle{empty}	

\noindent This version of the article has been accepted for publication, after a peer-review
process, and is subject to Springer Nature’s AM terms of use, but is not the Version
of Record and does not reflect post-acceptance improvements, or any corrections.
The Version of Record is available online at:

\noindent \url{https://doi.org/10.1007/s40314-020-01356-8}

\newpage

\setcounter{page}{1}

\maketitle

\begin{abstract}
An efficient approximate version of  implicit Taylor methods for initial-value problems of systems of
ordinary differential equations (ODEs) is introduced. The approach,
based on an approximate formulation of  Taylor methods, produces a
method that requires less evaluations of the function that defines the
ODE and its derivatives than the usual version. On the other hand, an
efficient numerical solution of the equation that arises from the
discretization by means of Newton's method is introduced for an
implicit scheme of any order. Numerical experiments illustrate that
the resulting algorithm is simpler to implement and has better
performance than its exact counterpart.


\end{abstract}

\section{Introduction}

\subsection{Scope} This work is related to numerical methods for the solution of the autonomous system of ordinary differential equations (ODEs)
\begin{align} \label{goveq}  \begin{split} 
 &  \boldsymbol{u}'(t)= \boldsymbol{f} \bigl( \boldsymbol{u}(t) \bigr), \quad t\in (t_0, T], \\
&    \boldsymbol{u} (t) = \bigl(u_1(t), \dots, u_M(t) \bigr)^{\mathrm{T}}, \quad 
 \boldsymbol{f} ( \boldsymbol{u}) =  \bigl(f_1( \boldsymbol{u}), \dots, f_M( \boldsymbol{u} ) \bigr)^{\mathrm{T}},  \end{split} 
\end{align} 
where derivatives of a vector of univariate scalar functions are understood com\-po\-nent-wise, 
posed along with initial data $\boldsymbol{u}(t_0) = \boldsymbol{u}_0$.

Taylor series methods for the numerical solution of initial value problems of ODEs compute approximations to the solution of the ODE for the next time instant using a Taylor polynomial of the unknown. The resulting methods are simple, since the expressions required for the iteration are exactly computable (i.e., with no error) from the equation, and the truncation error is governed by the error term of the Taylor formula, so that the order of accuracy of the global error of the method corresponds to the degree of the Taylor polynomial used. 
However, their implementation depends on the terms
  involved in the Taylor series, i.e., derivatives of the right-hand
  side whose computation requires intensive symbolic
  calculus and are specific to each individual problem. Moreover, the need for solving auxiliary nonlinear equations, especially within the implicit versions, makes them computationally expensive, especially as the order of accuracy required increases.

In this work we focus on implicit Taylor methods, obtained
  by computing the Taylor polynomials centered on a future time
  instant, and often used to solve problems where
  explicit methods have strong stability restrictions, in particular
  stiff systems of ODEs (Haier and Wanner, 1996). First of all we apply a
  strategy, based on the work by Baeza et al.\ (2017) 
  for the explicit Taylor method,
  to efficiently approximate the derivatives of~$\boldsymbol{f}$. This approximation inherits the ease of
  implementation and performance of the explicit version.

The implicit character of the method requires the
  solution of an auxiliary system of equations, usually by Newton's
  method, which requires the computation of the  Jacobian
  matrix. This may be an easy task for low-order methods, but the
  resulting iteration can become complicated as the order of the
  scheme increases. We propose a new formulation to obtain high-order
  implicit Taylor schemes that are simpler to implement and more
  efficient than the exact implicit Taylor methods, which compute
  derivatives symbolically. This is the main novelty of this work.

That said, we remark that it is not our purpose to
  present a numerical scheme that can compete with 
  	any implicit scheme in any situation, but to introduce a
  methodology to obtain $R$-th order implicit Taylor schemes for
  systems of $M$~scalar ODEs, with arbitrarily high $M\in \mathbb{N}$, that can be easily implemented
  and efficiently solved, independently of the complexity of the  function~$\boldsymbol{f}$, thus removing the leading difficulty of exact implicit Taylor methods. Very-high-order implicit methods are a must in some problems: for instance, in dynamical systems and mechanics, there is
  a need of high-order (at least, greater than 12) ODE solvers,
  especially Taylor integrators, as exposed for instance by Jorba and Zou (2005), 
   Barrio et al.\ (2011), and Abad et al.\ (2012).

\subsection{Related work}

Miletics and Moln\'{a}rka (2004)  propose an alternative based on a numerical approximation of the derivatives of $f$ in ODEs of the form $u'=f(u)$ for the explicit Taylor method up to fourth order and in 
Miletics and Moln\'{a}rka (2005)  for the implicit version up to fifth order. Later on, in 
 Baeza et al.\ (2017), a procedure to obtain a numerical approximation of $f(u)=f\circ u$ was presented to generate arbitrarily high order Taylor schemes, inspired by an approximate Cauchy-Kovalevskaya procedure developed for systems of conservation laws 
 by Zor\'{\i}o et al.\ (2017), which simplifies the exact version presented by Qiu and Shu (2003). The method presented 
  by Baeza et al.\ (2017)   relies on the approximate computation of the terms that appear in the Taylor polynomials, in terms of function evaluations only, avoiding the explicit computation of the derivatives, leading to a method which is simple to implement and outperforms its exact counterpart for complex systems.

Further references to implicit Taylor methods addressing combinations of implicit and explicit 
 steps to improve stability or accuracy include  Kirlinger and Corliss (1991)  and Scott (2000).

\subsection{Outline of the paper}  The work is organized as follows: In Section \ref{sec:taylor} the 
 basic facts about the exact Taylor methods are reviewed. A  general procedure to generate Taylor schemes of arbitrarily high accuracy order through Fa\`a di Bruno's formula  (Fa\`a di Bruno, 1855) is described, as well as its corresponding approximate version presented by Baeza et al.\ (2017).   
Section~\ref{sec:imptaylor} is devoted to the description of the novel
formulation of implicit Taylor methods, following an idea akin
to~Baeza et al.\ (2017). Section \ref{sec:Newtoniter}
describes an efficient implementation of the Newton iteration required
to update the solution of implicit Taylor methods. 
Section \ref{sec:numexp} stands for several numerical experiments in which the approximate version of the implicit Taylor methods is compared against its exact counterpart, as well as against the approximate explicit version. Finally, in Section \ref{sec:conclusions} some conclusions are drawn.

\section{Taylor methods}\label{sec:taylor}

\subsection{Preliminaries} 
The (explicit) $R$-th order Taylor 
methods are based on the expansion of the unknown function 
\begin{equation}\label{eq:poltay}
\boldsymbol{u} (t+h) = 
  \boldsymbol{u} (t)+h \boldsymbol{u}'(t)+\frac{h^2}{2} \boldsymbol{u}''(t)+\dots+\frac{h^R}{R!}\boldsymbol{u}^{(R)}(t) +\frac{h^{R+1}}{(R+1)!}\boldsymbol{u}^{(R+1)}(\xi)
\end{equation}
with $\xi$ belonging to the open interval $I(t, t+h)$ defined by~$t$ and~$t+h$. This expansion is  valid provided
$u_1, \dots, u_{M}  \in \smash{\mathcal{C}^{R}(\bar{I}(t, t+h))}$ and $\smash{u_1^{(R+1)}}, \dots, \smash{u_M^{(R+1)}}$ are bounded  in $I(t, t+h)$, where $\bar{I}(t, t+h)$
denotes the closure of $I(t, t+h)$. 
Consider an equally spaced set of $N+1$ points 
$t_n = t_0 + nh$, $0\leq n \leq N$, $h=T/{N}$.
Dropping the last term in \eqref{eq:poltay} and taking $t=t_n$ one obtains the approximation
\begin{align}\label{eq:taylor}
  \boldsymbol{u}(t_n+h)= \boldsymbol{u} (t_{n+1})\approx
  \boldsymbol{u}(t_n)+h \boldsymbol{u}'(t_n)+\frac{h^2}{2} \boldsymbol{u}''(t_n)+\dots+\frac{h^R}{R!} \boldsymbol{u}^{(R)}(t_n). 
\end{align}   
Then   \eqref{goveq} can be used to write  
\begin{align}  \label{uktn} 
\boldsymbol{u}^{(k)}(t_n)= \bigl( \boldsymbol{f}( \boldsymbol{u}) \bigr)^{(k-1)}(t_n)
= \frac{\mathrm{d}^{k-1}}{\mathrm{d} t^{k-1}}  \bigl( \boldsymbol{f}( \boldsymbol{u}(t)) \bigr)\biggr|_{t= t_n} 
, \quad  1\leq k\leq R. 
\end{align} 
Consequently, the first step to apply Taylor methods is to compute these derivatives up to an appropriate order.

\subsection{Fa\`{a} di Bruno's formula} \label{subsec:fdb}  The   evaluation 
of high-order derivatives of 
 the function $t \mapsto ( \boldsymbol{f} \circ \boldsymbol{u} ) (t)$, 
  which arise in \eqref{eq:taylor},   is greatly simplified by   Fa\`{a} di Bruno's formula, 
  as stated by  Baeza et al.\ (2017). To this end, we recall that 
   for a multi-index  $\boldsymbol{s} = (s_1, \dots, s_r) \in \mathbb{N}_0^r$,
    one defines $|\boldsymbol{s} | := s_1 + \dots + s_r$ and
    \begin{align*} 
     \binom{r}{\boldsymbol{s}} := \frac{r!}{s_1! s_2 ! \cdots s_r!}. 
     \end{align*} 
   Moreover, for $r \in \mathbb{N}$   we define an index set 
   \begin{align*} 
    \mathcal{P}_r := \left\{ \boldsymbol{s} \in \mathbb{N}_0^r \left|  \sum_{ \nu =1}^r \nu s_\nu  =r 
    \right. \right\},  
    \end{align*} 
    and $(\boldsymbol{D}^{\boldsymbol{s}} \boldsymbol{u})(t)$ to be a matrix 
     of size $M \times | \boldsymbol{s}|$ whose $(s_1 + \dots+s_{j-1}+i)$-th column
      is given by 
      \begin{align}  \label{dsdef} 
       \bigl( (\boldsymbol{D}^{\boldsymbol{s}} \boldsymbol{u})(t)\bigr)_{s_1 + \dots+s_{j-1}+i} 
        = \frac{1}{j!} \frac{\mathrm{d}^j}{\mathrm{d} t^j} \boldsymbol{u} (t), 
        \quad  i=1 , \dots, s_j, \quad j= 1, \dots, r.  
        \end{align}  
   Finally, we denote by $f^{(k)} \bullet \boldsymbol{A}$ the action of the $k$-th order derivative 
       tensor of~$f$ on an $M \times k$ matrix $\boldsymbol{A} = (A_{ij})$: 
        \begin{align*} 
         f^{(k)} \bullet \boldsymbol{A} 
          := \sum_{i_1, \dots, i_k=1}^{M} \frac{\partial^k f}{\partial u_{i_1} 
           \cdots \partial u_{i_k} } ( \boldsymbol{u} ) A_{i_1, 1} \cdots 
            A_{i_k, k}. 
            \end{align*}

\begin{proposition}[Fa\`{a} di Bruno's formula (Fa\`{a} di Bruno, 1855)] \label{prop:fdb}  Assume 
 that the functions $f:\mathbb{R}^M \to \mathbb{R}$ and $\boldsymbol{u}: 
 \mathbb{R} \to \mathbb{R}^M$ are $r$ times continuously differentiable. Then 
 \begin{align}\label{fdbformula} 
  \frac{\mathrm{d}^r}{\mathrm{d} t ^r} 
   f \bigl( \boldsymbol{u}(t) \bigr)  \equiv  \bigl(f(\boldsymbol{u}) \bigr)^{(r)} (t) 
    = \sum_{\boldsymbol{s} \in \mathcal{P}_r} \binom{r}{\boldsymbol{s}} 
     \Bigl( \bigl(f ( \boldsymbol{u}) \bigr)^{(|\boldsymbol{s}|)} 
      \bullet (\boldsymbol{D}^{\boldsymbol{s}} \boldsymbol{u})\Bigr) (t). 
\end{align} 
\end{proposition} 

  Proposition~\ref{prop:fdb} applies to just one scalar function~$f$, so to
 obtain  all components of, say, $\smash{(\boldsymbol{f} ( \boldsymbol{u}))^{(k)}} (t_n)$ 
  in \eqref{uktn}, we must apply \eqref{fdbformula} to each of the components of 
   $\boldsymbol{f} = (f_1, \dots, f_M)^{\mathrm{T}}$. Clearly, the matrix 
    $\boldsymbol{D}^{\boldsymbol{s}} \boldsymbol{u}$ is the same for all these components.

\subsection{Explicit Taylor methods}
The derivatives $\smash{( \boldsymbol{f}( \boldsymbol{u}) )^{(k-1)}}$ can be evaluated  by using Fa\`a di Bruno's formula \eqref{fdbformula}  (see Baeza et al.\ (2017) 
 for more details), leading to an expression
 of  $\smash{\boldsymbol{u}^{(k)}(t_n)}$ in terms  of $\boldsymbol{u} (t_n)$ and derivatives 
  of~$\boldsymbol{f}$, namely  
\begin{align}\label{eq:ck1} \begin{split} 
\boldsymbol{u}^{(k)}(t_n)  & = \boldsymbol{G}_k 
 \Bigl( \boldsymbol{u}(t_n), \bigl( \boldsymbol{f} (  \boldsymbol{u}) \bigr)     (t_n), 
 \bigl( \boldsymbol{f} (  \boldsymbol{u}) \bigr)'     (t_n),  \dots , 
  \bigl( \boldsymbol{f} (  \boldsymbol{u}) \bigr)^{(k-1)}      (t_n)
  \Bigr)=\boldsymbol{\tilde {G}}_k \bigl(\boldsymbol{u}(t_n) \bigr).  \end{split} 
\end{align}
Replacing the derivatives $\smash{\boldsymbol{u}^{(k)} (t_n)}$ in \eqref{eq:taylor} 
 by \eqref{eq:ck1}, we obtain the  expression
\begin{align}\label{eq:tr}
  \boldsymbol{u} (t_{n+1}) \approx T_R \bigl( \boldsymbol{u}(t_n), h \bigr) 
   = \boldsymbol{u} (t_n)+\sum_{k=1}^{R}
    \frac{h^k}{k!} \boldsymbol{\tilde{G}}_k \bigl( \boldsymbol{u}(t_n) \bigr). 
\end{align}

The  $R$-th order Taylor method
\begin{align}  \label{ROT} 
 \boldsymbol{u}_{n+1} = 
  T_R( \boldsymbol{u}_n, h)
  \end{align} 
  is  then obtained  by replacing  the exact values of  the solution $u(t_n)$ and $u(t_{n+1})$ by their corresponding approximations in \eqref{eq:tr},  
denoted by~$\boldsymbol{u}_n$ and~$\boldsymbol{u}_{n+1}$, respectively. This means that 
 the following expression is used in \eqref{ROT}: 
\begin{align}\label{eq:pep1}
    T_R( \boldsymbol{u}_n, h) = \boldsymbol{u}_n+\sum_{k=1}^{R}
    \frac{h^k}{k!} \boldsymbol{u}_n^{(k)}, \quad \boldsymbol{u}_n^{(k)}:= 
    \boldsymbol{\tilde{G}}_k(\boldsymbol{u}_n).
\end{align}
From \eqref{eq:poltay}  and \eqref{eq:pep1} we infer that  the local truncation error is given by 
$$
\boldsymbol{E}_{\mathrm{L}} = \frac{h^{R+1}}{(R+1)!} \boldsymbol{u}^{(R+1)}(\xi),  
$$
so that     $\boldsymbol{E}_{\mathrm{L}} = \mathcal{O}(h^{R+1})$ as long as $\boldsymbol{u}^{(R+1)}$ is bounded in $[t_0, T]$. One then obtains that the method \eqref{ROT} has an $\mathcal{O}(h^{R})$   global error .

\section{Implicit Taylor methods}\label{sec:imptaylor}

\subsection{Exact implicit Taylor methods} 
Implicit Taylor methods are based on approximating 
$\boldsymbol{u}(t_n)$ by means of the Taylor polynomial of~$\boldsymbol{u}$ centered at $t_{n+1}$: 
\begin{align} \label{eq:taylor_implicit-prel}
  \boldsymbol{u}(t_n) \approx T_R \bigl( \boldsymbol{u}(t_{n+1}), -h \bigr),      \end{align}
so that the value of  $ \boldsymbol{u}_{n+1} \approx  \boldsymbol{u}(t_{n+1})$ 
 is determined as  solution of the nonlinear system of algebraic equations 
\begin{align}\label{eq:taylor_implicit}
  \boldsymbol{u}_n = T_R( \boldsymbol{u}_{n+1}, -h).
\end{align} 
In the easiest case, with $R=1$, one gets the implicit Euler method. As in the case of explicit Taylor methods, 
the expressions of $\smash{\boldsymbol{u}^{(k)}(t_{n+1})}$ that appear in \eqref{eq:taylor_implicit-prel} can be expressed as  functions 
of $\boldsymbol{u}(t_{n+1}$) and the derivatives of~$\boldsymbol{f}$. As an example, the second-order implicit Taylor method
is given by 
\begin{align}\label{eq:tay_implicit_2}
\boldsymbol{u}_n = \boldsymbol{u}_{n+1} - h \boldsymbol{f}(\boldsymbol{u}_{n+1}) + \frac{h^2}{2}\left(\frac{\partial \boldsymbol{f}}{\partial \boldsymbol{u}}( \boldsymbol{u}_{n+1}) 
 \boldsymbol{f}(\boldsymbol{u}_{n+1})\right), 
\end{align}
 where $\smash{\partial \boldsymbol{f} / \partial \boldsymbol{u}  = (\partial f_i / \partial u_j)_{1 \leq i,j \leq  M}}$ is the Jacobian matrix of $\boldsymbol{f} ( \boldsymbol{u})$.  In what follows, the family of methods 
   based on \eqref{eq:taylor_implicit} will be referred to as {\em exact} implicit Taylor methods 
    since they are based on exact expressions of the derivatives of~$\boldsymbol{f}$.

\subsection{Approximate implicit Taylor methods}
  Let us briefly review approximate explicit Taylor methods 
 as described by Baeza et al.\ (2017),
   whose formulation will be used to motivate and introduce our novel  \textit{approximate} implicit
  Taylor (henceforth, AIT) methods. These methods 
  are based on computing  approximations of the derivatives 
in \eqref{eq:taylor} by means of finite differences, so that $\smash{\boldsymbol{u}^{(k)}(t_n)}$ is replaced 
by an approximation 
\begin{align*} 
\boldsymbol{v}_{h,n}^{(k)} = \boldsymbol{u}^{(k)}(t_n) +\mathcal{O}(h^{R-k+1}),  
 \quad k=2,\dots, R, 
 \end{align*} 
resulting in an $R$-th order accurate  method 
\begin{align*} 
\boldsymbol{v}_{h, n+1}= \boldsymbol{v}_{h,n} + \sum_{k=1}^R   \frac{h^k}{k!}     \boldsymbol{v}_{h,n}^{(k)},
\end{align*} 
where the approximations $\smash{v_{h,n}^{(k)}}$ are computed 
as follows:
\begin{align*}
\boldsymbol{v}_{h,n}^{(0)} &= \boldsymbol{u}_n,\\
\boldsymbol{v}_{h,n}^{(1)} &= \boldsymbol{f}(\boldsymbol{u}_n),\\
\boldsymbol{v}_{h,n}^{(k+1)} &= \Delta_{h}^{k,\lceil\frac{R-k}{2}\rceil} \boldsymbol{f}
 \bigl( \boldsymbol{P}_n^k(h) \bigr), \quad k= 1, \dots, R-1. 
\end{align*}
Here we recall  that $\lceil \cdot \rceil$ denotes the so-called {\em ceiling} operator 
defined by $\lceil x
  \rceil=\min \{n \in \mathbb{Z}  \mid  x\leq n\}$. Moreover,  
$\boldsymbol{P}^k(\rho)$
	 is the  $M$-component  vector given by 
\begin{align*} 
P^k_n(\rho) = \sum_{l=0}^k\frac{\boldsymbol{v}^{(l)}_{h,n}}{l!}\rho^l, \quad n=1, \dots, M, 
\end{align*} 
and 
$\Delta_{h}^{p,q}$ is the centered finite-difference operator that approximates $p$-th order derivatives to order~$2q$
on a grid with spacing $h$, i.e., the one that satisfies  
$$\Delta_h^{p,q}(y) = y^{(p)} + \mathcal{O}(h^{2q})$$
for a sufficiently differentiable function~$y$.  (The operator 
$\Delta_{h}^{p,q}$ is understood as acting on each component of 
$\boldsymbol{f}
 ( \boldsymbol{P}_n^k(h))$ separately.)

 There exist constants $\smash{\beta^{k, R}_j}$ so that for some integers $\gamma_{k,R}$, we can write (see  Zor\'{\i}o et al., 2017)   
\begin{equation}\label{eq:vk}
\boldsymbol{v}_{h,n}^{(k+1)}=h^{-k}\sum_{j=-\gamma_{k,R}}^{\gamma_{k,R}}
  \beta_j^{k,R} \boldsymbol{f} \left(\sum_{l=0}^k \frac{(jh)^l}{l!} \boldsymbol{v}_{h,n}^{(l)}\right). 
  \end{equation}
Using these approximations of the derivatives, and with the notation of the previous sections, one obtains the approximate  explicit Taylor method 
\begin{align}\label{eq:tay_apr}
  \boldsymbol{u}_{n+1} = \tilde{T}_R( \boldsymbol{u}_{n}, h).
\end{align}
For instance, the second-order approximate Taylor method is based on the
approximation 
\begin{align*} \boldsymbol{u}^{(2)}(t_n)= \bigl( \boldsymbol{f}( \boldsymbol{u}) \bigr)'(t_n)\approx
\frac{1}{2h} \Bigl( \boldsymbol{f} \bigl( \boldsymbol{u}(t_n)+h \boldsymbol{f} ( \boldsymbol{u}(t_n))\bigr)- \boldsymbol{f} \bigl( \boldsymbol{u}(t_n)-h \boldsymbol{f}(u(t_n))\bigr) \Bigr),
\end{align*} 
hence the method can be written as
\begin{align*}
  \boldsymbol{u}_{n+1}=\boldsymbol{u}_n+h \boldsymbol{f}( \boldsymbol{u}_n)+\frac{h}{4}
   \Big(\boldsymbol{f}( \boldsymbol{u}_n+h \boldsymbol{f}(u_n)\bigr)-
    \boldsymbol{f} \bigl(\boldsymbol{u}_n-h \boldsymbol{f}(u_n)\bigr)\Bigr), 
\end{align*}  
i.e.,
\begin{align*}
  \tilde{T}_2(\boldsymbol{u}_n, h) =\boldsymbol{u}_n+h \boldsymbol{f}( \boldsymbol{u}_n)+\frac{h}{4}
   \Big(\boldsymbol{f}( \boldsymbol{u}_n+h \boldsymbol{f}(u_n)\bigr)-
    \boldsymbol{f} \bigl(\boldsymbol{u}_n-h \boldsymbol{f}(u_n)\bigr)\Bigr). 
\end{align*}
The new methods advanced in this contribution, namely 
{\em approximate implicit} Taylor methods,  are obtained by  replacing~$h$ by~$-h$  and 
 interchanging~$\boldsymbol{u}_n$ and~$\boldsymbol{u}_{n+1}$  in 
\eqref{eq:tay_apr}:
\begin{align*}
  \boldsymbol{u}_{n} =\tilde{T}_R( \boldsymbol{u}_{n+1}, -h).
\end{align*}
For the case of second order ($R=2$), the implicit second-order approximate Taylor method is 
\begin{align}\label{eq:tay_apr_implicit_2} 
 \boldsymbol{u}_n&=\tilde{T}_2( \boldsymbol{u}_{n+1}, -h) \\ & = \boldsymbol{u}_{n+1}-hf( \boldsymbol{u}_{n+1}) -\frac{h}{4} \Bigl( \boldsymbol{f}
 \bigl( \boldsymbol{u}_{n+1}-h  \boldsymbol{f}(  \boldsymbol{u}_{n+1}) \bigr)-  \boldsymbol{f} \bigl(  \boldsymbol{u}_{n+1}+h  \boldsymbol{f}(  \boldsymbol{u}_{n+1}) \bigr) \Bigr). \nonumber  
\end{align}

\subsection{Linear stability} The  linear stability of a numerical
scheme for initial value problems of ordinary differential equations
is usually examined by applying it to the   scalar linear equation 
\begin{align} \label{lineq} 
  u'=\lambda u, \quad  \lambda \in \mathbb{C}, \quad \mathrm{Re} \, \lambda< 0.
\end{align}
For the sake of completeness,
we consider the non-homogeneous linear ODE
\begin{equation*}
  u'=\lambda u+g(t), \quad  \lambda \in \mathbb{C},
\end{equation*}
with $g$ sufficiently  smooth.
For the solution $u$ of the ODE, we can establish by induction on $k$
that 
\begin{equation*}
  u^{(k)}=\lambda^k u+\sum_{j=0}^{k-1}\lambda^{k-j-1}g^{(j)}(t),
\end{equation*}
so the explicit Taylor method reads in this case as 
\begin{align*}  u_{n+1}&=\sum_{k=0}^{R} \frac{h^{k}}{k!}
  \Bigg(\lambda^k
  u_{n}+\sum_{j=0}^{k-1}\lambda^{k-j-1}g^{(j)}(t_{n})\Bigg)\\
    &=u_{n}\sum_{k=0}^{R} \frac{(h\lambda)^{k}}{k!} +
  \sum_{j=0}^{R-1} \frac{g^{(j)}(t_{n})}{\lambda^{j+1}}
  \sum_{k=j+1}^{R}\frac{(h\lambda)^{k}}{k!}\\
  &=Q_{R}(h\lambda)u_{n} +
  \sum_{j=0}^{R-1} \frac{g^{(j)}(t_{n})}{\lambda^{j+1}}
  \big(Q_{R}(h\lambda)-Q_{j}(h\lambda)\big),
\end{align*}  
where 
\begin{equation*}
  Q_{j}(x)=\sum_{k=0}^{j}\frac{x^k}{k!}.
\end{equation*}

The implicit Taylor method is obtained by interchanging the roles of $n$
and $n+1$ and  reads as 
\begin{align}
\notag
  u_{n}&=Q_{R}(-h\lambda)u_{n+1} +
  \sum_{j=0}^{R-1} \frac{g^{(j)}(t_{n+1})}{\lambda^{j+1}}
  \big(Q_{R}(-h\lambda)-Q_{j}(-h\lambda)\big),\\
\label{eq:pep1234}
  u_{n+1}&=\frac{1}{Q_R(-h\lambda)} u_{n}-\sum_{j=0}^{R-1}
  \frac{g^{(j)}(t_{n+1})}{\lambda^{j+1}} \biggl(1-\frac{Q_j(-h\lambda)}{Q_R(-h\lambda)} \biggr).
\end{align}
In particular, for $g=0$, the explicit and implicit Taylor methods of
order $R$ are given by the respective expressions
\begin{equation} \label{expltayl} 
   u_{n+1}=Q_R(h\lambda) u_n
 \end{equation}
 and   
   \begin{equation*}   u_{n+1}=\frac{1}{Q_R(-h\lambda)} u_n.
 \end{equation*}

The exact Taylor method of order $R$ is stable provided that 
$|Q_R(h\lambda)| < 1$. Since $\mathrm{Re} \, \lambda <0$, this condition is usually 
satisfied on a bounded domain only (as can be inferred from $R=1$, in
which case 
 \eqref{expltayl} is the explicit Euler  method). On the other hand, 
  the exact implicit Taylor method
is stable for those values of $z = h \lambda$  that satisfy 
\begin{align*} \begin{split} 
z \in \mathcal{S} := & \bigl\{ z \in \mathbb{C} \mid \mathrm{Re} \, z < 0,  |Q_R(-z)|^{-1} < 1 \bigr\}=\bigl\{ z \in \mathbb{C} \mid \mathrm{Re} \, z < 0,  |Q_R(-z)| > 1 \bigr\}. \end{split} 
\end{align*}

 As for its exact counterpart, in  Baeza et al.\ (2020)  it is shown that 
the approximate explicit Taylor method applied to \eqref{lineq}  is
$\tilde{T}_R(u_{n}, h)=Q_R(h\lambda) u_n$, thus the implicit version is $\tilde{T}_R (u_{n}, -h)=Q_R(-h\lambda)^{-1} u_n$,
and therefore both methods have the same stability region 
as their corresponding exact versions, in particular
the approximate implicit Taylor method is absolutely stable 
whenever $\lambda < 0$.

\section{Newton iteration}\label{sec:Newtoniter}

The computation of $\boldsymbol{u}_{n+1}$ for given  $\boldsymbol{u}_n$ using an implicit 
 method requires the solution of an
auxiliary equation  
$\boldsymbol{F}( \boldsymbol{u}_{n+1})= \boldsymbol{0}$, 
which is often approximated by Newton's method. 
In this section, we address the computation of the elements required for 
Newton's method for both the exact and approximate implicit Taylor methods, which will lead to a new, more efficient formulation for the approximate scheme. Although line-search strategies   (Dennis and Schnabel, 1996)  for damping Newton iteration can be used to enhance global convergence we have not used them in our experiments.

\subsection{Exact implicit Taylor method}
As an example, let us consider the  scalar nonlinear problem
\begin{equation*}
  u'=u+u^2 \Rightarrow u''=(1+2u)u'=(1+2u)(u+u^2).
\end{equation*}
The second-order exact implicit Taylor method can be written as
\begin{equation*}
  u_n=u_{n+1}-h(u_{n+1}+u_{n+1}^2)+\frac{h^2}{2}(1+2u_{n+1})(u_{n+1}+u_{n+1}^2),
\end{equation*}
which  requires the solution of the 
 following cubic equation:
\begin{equation}\label{eq:newton1}
  F(u_{n+1}):=u_{n+1}-h(u_{n+1}+u_{n+1}^2)+\frac{h^2}{2}(1+2u_{n+1})(u_{n+1}+u_{n+1}^2)-u_n=0.
\end{equation}
In the general case, the solution of $\boldsymbol{F}( \boldsymbol{u}_{n+1})= \boldsymbol{0}$   by means of Newton's method requires the computation of the derivative $F'(u_{n+1})$. 
In the case of \eqref{eq:newton1} this is an easy task, but in general the 
resulting iteration can become complicated. 

To simplify the computation of the Jacobian matrix, we introduce
\begin{align*} 
\boldsymbol{z}_k\approx \boldsymbol{u}_{n+1}^{(k)}, \quad  k=0,\dots,R, 
\end{align*} 
and use Fa\`a di Bruno's formula \eqref{fdbformula}  to get the system
\begin{align*}
  \boldsymbol{u}_{n}&= \boldsymbol{z}_0-h \boldsymbol{z}_1+\dots+(-1)^R\frac{h^R}{R!}\boldsymbol{z}_R,\\
  \boldsymbol{z}_1&=\boldsymbol{f}(\boldsymbol{z}_0),\\
  \boldsymbol{z}_{n+1}&=\sum_{ \boldsymbol{s} \in\mathcal{P}_r}
    \binom{r}{\boldsymbol{s}} \begin{pmatrix} 
     f_1^{(|\boldsymbol{s}|)}(\boldsymbol{z}_0) \bullet \boldsymbol{\tilde{D}}^{\boldsymbol{s}}
      \boldsymbol{z} \\  \vdots \\ f_M^{(|\boldsymbol{s}|)}(\boldsymbol{z}_0) \bullet 
      \boldsymbol{\tilde{D}}^{\boldsymbol{s}}
      \boldsymbol{z}  \end{pmatrix} ,\quad   r=1,\dots,R-1,
  \end{align*}  
     where  the definition of $\smash{\boldsymbol{\tilde{D}}^{\boldsymbol{s}}
      \boldsymbol{z}}$  mimics that of $\smash{\boldsymbol{D}^{\boldsymbol{s}}
      \boldsymbol{z}}$ in \eqref{dsdef}, by taking into account that 
      $\boldsymbol{z}_k \approx  
      \smash{ \boldsymbol{u}^{(k)}} (t)$, namely 
    \begin{align*}       \bigl( \boldsymbol{\tilde{D}}^{\boldsymbol{s}} \boldsymbol{z} \bigr)_{s_1 + \dots+s_{j-1}+i} 
        = \frac{1}{j!}  \boldsymbol{z}_j,  
        \quad  i=1 , \dots, s_j, \quad j= 1, \dots, r.  
        \end{align*}
   These equations   can be differentiated systematically. For instance, for the case of one 
   scalar equation, $M=1$, one gets 
\begin{align*}  
& \partial_{z_0} 
 \Biggl( \sum_{\boldsymbol{s} \in \mathcal{P}_r} \binom{r}{\boldsymbol{s}}  
  f^{(|\boldsymbol{s}|)} (z_0) \boldsymbol{D}^{\boldsymbol{s}}  z \Biggr)  = \sum_{\boldsymbol{s}  \in \mathcal{P}_r} 
   \binom{r}{\boldsymbol{s}}  f^{(|\boldsymbol{s}|+1)} (z_0) \Bigl( \frac{z_1}{1!} \Bigr)^{s_1} \cdots 
    \Bigl( \frac{z_r}{r!} \Bigr)^{s_r}, \\
 &    \partial_{z_j} 
 \Biggl( \sum_{\boldsymbol{s} \in \mathcal{P}_r} \binom{r}{\boldsymbol{s}}  
  f^{(|\boldsymbol{s}|)} (z_0) \boldsymbol{D}^{\boldsymbol{s}}  z \Biggr) \\
  & \qquad \qquad = \sum_{\boldsymbol{s}  \in \mathcal{P}_r}  \frac{s_j}{j!} 
   \binom{r}{\boldsymbol{s}}  f^{(|\boldsymbol{s}|+1)} (z_0) \Bigl( \frac{z_1}{1!} \Bigr)^{s_1} \cdots 
     \Bigl( \frac{z_j}{j!} \Bigr)^{s_j-1}  \cdots 
    \Bigl( \frac{z_r}{r!} \Bigr)^{s_r}. 
\end{align*} 
For this scalar  case and the second-order implicit Taylor method \eqref{eq:tay_implicit_2}
the system to be solved is
\begin{align} \label{syst:ex1}  \begin{split} 
0&=z_0 -hz_1+\frac{h^2}{2}z_2 - u_n,\\
0&=f(z_0)-z_1,\\
0 &= f'(z_0) z_1 -z_2. 
\end{split}
\end{align}  
If we rewrite system \eqref{syst:ex1} as $\boldsymbol{F} ( z_0, z_1, z_2) = \boldsymbol{0}$, with the function $\boldsymbol{F}$ defined by
\begin{align*} 
\boldsymbol{F}( z_0, z_1, z_2)= 
\begin{pmatrix} 
F_1 ( z_0, z_1, z_2)\\
F_2 ( z_0, z_1, z_2)\\
F_3 ( z_0, z_1, z_2)
\end{pmatrix}&=
\begin{pmatrix} 
\displaystyle z_0 -hz_1+\frac{h^2}{2}z_2 - u_n,\\
 f(z_0)-z_1,\\
 f'(z_0) z_1 -z_2
\end{pmatrix},
\end{align*}   
then the corresponding Jacobian matrix is
\begin{equation} \label{eq:jacob2}
\displaystyle \mathcal{J}_{\boldsymbol{F}} ( z_0, z_1, z_2) = 
 \begin{bmatrix}
\displaystyle\frac{\partial F_1}{\partial z_0}& \displaystyle\frac{\partial F_1}{\partial z_1} & \displaystyle\frac{\partial F_1}{\partial z_2} \\[3mm]
\displaystyle\frac{\partial F_2}{\partial z_0} & \displaystyle\frac{\partial F_2}{\partial z_1} &\displaystyle \frac{\partial F_2}{\partial z_2}\\[3mm]
\displaystyle\frac{\partial F_3}{\partial z_0} & \displaystyle\frac{\partial F_3}{\partial z_1} & \displaystyle\frac{\partial F_3}{\partial z_2}\\
\end{bmatrix} =
  \begin{bmatrix}
    1&-h & h^2 / 2 \\
    f'(z_0) & -1 & 0\\
    f''(z_0)z_1 & f'(z_0) &-1\\
\end{bmatrix}.
\end{equation}

Depending on the expression of $f$ the Jacobian matrix may become 
highly complicate, even for low values of $R$.  It is clear that for
higher order methods, the system to be solved will be more
complicated. For instance, for $R=4$ it reads as 
\begin{align*} \begin{split} 
0&=z_0 -hz_1+\frac{h^2}{2}z_2 - \frac{h^3}{6}z_3 + \frac{h^4}{24} z_4 -u_n,\\
0&=f(z_0)-z_1,\\
0 &= f'(z_0) z_1 -z_2,\\
0&= f''(z_0)z_1^2+f'(z_0)z_2-z_3,\\
0&= f'''(z_0)z_1^3+3f''(z_0) z_1 z_2 +z_3 f'(z_0)-z_4,
\end{split}
\end{align*}  
which results in the expression
\begin{align}\label{eq:jacob4} \begin{split} 
&  \mathcal{J}_{\boldsymbol{F}} ( z_0, \ldots, z_4) \\ 
& = 
\begin{bmatrix} 
1                & -h         & \displaystyle{\frac{h^2}{2}}  & -\displaystyle{\frac{h^3}{6}} & 
 \displaystyle{\frac{h^4}{24}}  \\[4mm]  
 f'(z_0)     & -1         & 0                    & 0                  & 0 \\[2mm] 
f''(z_0)z_1 & f'(z_0) & -1                   & 0                   & 0 \\[2mm] 
f'''(z_0)z_1^2 +f''(z_0)z_2  & 2f''(z_0)z_1 & f'(z_0) & -1 & 0\\[2mm] 
f^{(4)}(z_0)z_1^3 +3 f'''(z_0)z_1z_2  & \multirow{2}{*}{$3 f'''(z_0) z_1^2 +3f''(z_0) z_2$}  &  \multirow{2}{*}{$3f''(z_0) z_1$} &  \multirow{2}{*}{$f'(z_0)$} &  \multirow{2}{*}{$-1$}\\ 
 +f''(z_0)z_3  &  & & & 
\end{bmatrix}. \end{split} 
\end{align}	
Note that the submatrix composed by the first three rows 
and columns of \eqref{eq:jacob4} is exactly \eqref{eq:jacob2}. It is easy to check that the 
Jacobian matrix corresponding to $R=3$ is the submatrix of \eqref{eq:jacob4} composed
by its first four rows and columns.

\subsection{Approximate implicit Taylor method}

For simplicity, let us start with the second-order approximate implicit Taylor method
\eqref{eq:tay_apr_implicit_2} for the scalar case $M=1$. Similarly to the exact case, we introduce the unknowns 
  $z_0 =u_{n+1}$, 
  $z_1 =f(u_{n+1})$ and 
\begin{align*}   
  z_2 =\frac12 \Bigl(f \bigl(u_{n+1}-hf(u_{n+1})\bigr)-f \bigl(u_{n+1}+hf(u_{n+1}) \bigr) \Bigr),
\end{align*}  
and one gets the system of equations
\begin{align*}
0&=z_0-hz_1-\frac{h}{2}z_2-u_n,\\
0&=f(z_0)-z_1,\\
0&=\frac12 f(z_0-hz_1)-\frac12 f(z_0+hz_1)-z_2,
\end{align*}
so that its solution  gives the terms that appear in \eqref{eq:tay_apr_implicit_2}. 

If we rewrite this system as $\boldsymbol{F} (z_0, z_1, z_2) = \boldsymbol{0}$ with the function $\boldsymbol{F}$ defined by
\begin{align*} 
\boldsymbol{F}( z_0, z_1, z_2)= 
\begin{pmatrix}
F_1 ( z_0, z_1, z_2)\\
F_2 ( z_0, z_1, z_2)\\
F_3 ( z_0, z_1, z_2)
\end{pmatrix}&=
\begin{pmatrix} 
\displaystyle z_0-hz_1-\frac{h}{2}z_2-u_n\\
f(z_0)-z_1\\
\frac12 f(z_0-hz_1)-\frac12 f(z_0+hz_1)-z_2
\end{pmatrix},
\end{align*}   
then the corresponding Jacobian  matrix (which is required so that  Newton's method  can be applied to this system) is now given by
\begin{align*} 
& \mathcal{J}_{\boldsymbol{F}} (z_0, z_1, z_2)   =  
 \begin{bmatrix}
\displaystyle\frac{\partial F_1}{\partial z_0}& \displaystyle\frac{\partial F_1}{\partial z_1} & \displaystyle\frac{\partial F_1}{\partial z_2} \\[3mm]
\displaystyle\frac{\partial F_2}{\partial z_0} & \displaystyle\frac{\partial F_2}{\partial z_1} &\displaystyle \frac{\partial F_2}{\partial z_2}\\[3mm]
\displaystyle\frac{\partial F_3}{\partial z_0} & \displaystyle\frac{\partial F_3}{\partial z_1} & \displaystyle\frac{\partial F_3}{\partial z_2} 
\end{bmatrix}\\ & = 
 \begin{bmatrix}
     1 & -h  &-  \displaystyle \frac{h}{2} \\[4mm] 
     f'(z_0)  &  -1  &  0  \\[2mm] 
   \displaystyle{\frac12 \bigl( f'(z_0-h z_1)- f'(z_0+h z_1) \bigr) } &
\displaystyle{-\frac{h}{2} \bigl( f'(z_0-h z_1)+f'(z_0+h z_1) \bigr)} & -1\\
  \end{bmatrix}.
\end{align*}

\subsection{General number of scalar equations} \label{subsec_newton}

Let us now consider the general case of a system of $M$~scalar 
ordinary differential equations. From \eqref{eq:vk}, the approximate implicit  $R$-th order Taylor method 
can be written as
\begin{align}\label{eq:1}
  \boldsymbol{u}_{n}&=\tilde{T}_R(\boldsymbol{u}_{n+1}, -h)=\sum_{k=0}^R\frac{(-h)^k}{k!}
  \boldsymbol{v}_{-h,n}^{(k)},\\
  \label{eq:2}
  \boldsymbol{v}_{-h,n}^{(k+1)}&=(-h)^{-k}\sum_{j=-\gamma_{k,R}}^{\gamma_{k,R}}
  \beta_j^{k,R} \boldsymbol{f} \left(\sum_{l=0}^k \frac{j^l (-h)^l}{l!} \boldsymbol{v}_{-h,n}^{(l)}\right).
\end{align}
Let us denote $ \boldsymbol{z}_k=(-h)^{k-1} \boldsymbol{v}_{-h,n}^{k}$,  
so that \eqref{eq:2} for $k-1$ reads as 
\begin{align*}
   \boldsymbol{z}_{k}&=\sum_{j=-\gamma_{k-1,R}}^{\gamma_{k-1,R}}
  \beta_j^{k-1,R}  \boldsymbol{f} \left(z_0-h\sum_{l=1}^{k-1} \frac{j^l}{l!} \boldsymbol{z}_l\right)
\end{align*}
and \eqref{eq:1} as
\begin{align*}
  \boldsymbol{u}_{n}&=\boldsymbol{z}_0-h\sum_{k=1}^{R}\frac{1}{k!} \boldsymbol{z}_k.
\end{align*}
We define the function  
$\smash{\boldsymbol{F} = ( \boldsymbol{F}_0, \boldsymbol{F}_1, \dots, 
 \boldsymbol{F}_M)^{\mathrm{T}}:  \quad \mathbb{R}^{(R+1)M}\to\mathbb{R}^{(R+1)M}}$   
  by
\begin{align*}
  \boldsymbol{F}_0&:=\boldsymbol{z}_0-h\sum_{k=1}^{R}\frac{1}{k!} \boldsymbol{z}_k- \boldsymbol{u}_n,\\
  \boldsymbol{F}_{k}&:=\sum_{j=-\gamma_{k-1,R}}^{\gamma_{k-1,R}}
  \beta_j^{k-1,R} \boldsymbol{f} \left(\boldsymbol{z}_0-h\sum_{l=1}^{k-1}  \frac{j^l}{l!} \boldsymbol{z}_l \right)-\boldsymbol{z}_k,\quad
  k=1,\dots,R.
\end{align*}
To solve $\boldsymbol{F}(\boldsymbol{z})=\boldsymbol{0}$ by Newton's method,  we compute the Jacobian matrix of~$\boldsymbol{F}$ as the block matrix 
\begin{align*} 
 \mathcal{J}_{\boldsymbol{F}} ( \boldsymbol{z}) = \bigl( \boldsymbol{F}_{i,j} (\boldsymbol{z}) \bigr)_{0 \leq i,j \leq R}, \quad \text{where} \quad   \boldsymbol{F}_{i,j} ( \boldsymbol{z} ) = 
  \frac{\partial \boldsymbol{F}_i}{\partial \boldsymbol{z}_j} ( \boldsymbol{z}) \in 
   \mathbb{R}^{M\times M}. 
 \end{align*} 
 If $\boldsymbol{I}_{M}$ denotes the 
  $M\times M$ identity matrix, we get
\begin{align*}
   \boldsymbol{F}_{0,0}&= \boldsymbol{I}_M,\\
  \boldsymbol{F}_{0,l} &=-\frac{h}{l!} \boldsymbol{I}_M,\quad l=1,\dots,R,\\
  \boldsymbol{F}_{k,0} &=\sum_{j=-\gamma_{k-1,R}}^{\gamma_{k-1,R}}\beta_j^{k-1,R} \boldsymbol{f}'\left(\boldsymbol{z}_0-h\sum_{l=1}^{k-1}\frac{j^l}{l!}
  \boldsymbol{z}_l\right),\quad  k=1,\dots,R,\\
   \boldsymbol{F}_{k,l}&=-h\sum_{j=-\gamma_{k-1,R}}^{\gamma_{k-1,R}} \beta_j^{k-1,R} \boldsymbol{f}'\left(\boldsymbol{z}_0-h\sum_{m=1}^{k-1}\frac{j^m}{m!}
  \boldsymbol{z}_m \right)\frac{j^l}{l!} ,\quad \begin{cases} l=1,\dots,k-1, \\ k=1,\dots,R, \end{cases} \\
  \boldsymbol{F}_{k,k}&=-\boldsymbol{I}_M,\quad k=1,\dots,R,\\
  \boldsymbol{F}_{k,l}&= \boldsymbol{0}, \quad l=k+1\dots, R, \quad  k=1,\dots,R.
\end{align*}  
Setting $\smash{\boldsymbol{\delta}^{(\nu)} = \boldsymbol{z}^{(\nu+1)} - \boldsymbol{z}^{(\nu)}}$,  we may write  an iteration of Newton's method as 
\begin{equation*}
 \mathcal{J}_{\boldsymbol{F}}  ( \boldsymbol{z}^{(\nu)}) \boldsymbol{\delta}^{(\nu)}=-
  \boldsymbol{F}( \boldsymbol{z}^{(\nu)}).
\end{equation*}  
In block form and dropping $\nu$, we get 
\begin{align*}
  \begin{bmatrix}
   \boldsymbol{F}_{0,0}  & \boldsymbol{F}_{0,1} & \cdots &\boldsymbol{F}_{0,R}  \\
    \boldsymbol{F}_{1,0}  &  \boldsymbol{F}_{1,1}  & \cdots & \boldsymbol{F}_{1,R} \\ 
    \vdots & \vdots & & \vdots \\ 
     \boldsymbol{F}_{R,0}  &  \boldsymbol{F}_{R,1}  & \cdots &\boldsymbol{F}_{R,R} 
  \end{bmatrix}
  \begin{pmatrix}
    \boldsymbol{\delta}_0\\
    \boldsymbol{\delta}_{1} \\ \vdots \\ \boldsymbol{\delta}_R 
  \end{pmatrix}
  &=
  -\begin{pmatrix}
   \boldsymbol{F}_0\\ \boldsymbol{F}_1\\ \vdots \\  \boldsymbol{F}_R 
  \end{pmatrix},
\end{align*}
which we write in compact form as 
\begin{align} \label{comp-syst} 
\begin{bmatrix} 
\boldsymbol{F}_{0,0} & \boldsymbol{F}_{0,1:R} \\ 
\boldsymbol{F}_{1:R,0} & \boldsymbol{F}_{1:R,1:R} \\
\end{bmatrix} \begin{pmatrix} 
 	\boldsymbol{\delta}_0 \\ \boldsymbol{\delta}_{1:R} \end{pmatrix} &=- 
	 \begin{pmatrix} 	\boldsymbol{F}_0 \\
 \boldsymbol{F}_{1:R} \end{pmatrix}.
\end{align}
Since $\boldsymbol{F}_{1:R,1:R}$ is blockwise lower triangular with the diagonal blocks 
 given by~$-\boldsymbol{I}_M$,  this matrix is invertible and 
we deduce that 
\begin{align*}
\boldsymbol{\delta}_{1:R}&=-  \boldsymbol{F}_{1:R, 1:R} ^{-1}( \boldsymbol{F}_{1:R}+
 \boldsymbol{F}_{1:R,0} \boldsymbol{\delta}_0),
\end{align*}
which, when inserted into the first equation of \eqref{comp-syst},  yields
\begin{align*}
  \boldsymbol{\delta}_0 
   = - \bigl( \boldsymbol{F}_{0,0} -   \boldsymbol{F}_{0,1:R} \boldsymbol{F}_{1:R, 1:R}^{-1} 
    \boldsymbol{F}_{1:R,0} \bigr)^{-1}  \bigl( 
     \boldsymbol{F}_{0} -   \boldsymbol{F}_{0,1:R} \boldsymbol{F}_{1:R, 1:R}^{-1} 
    \boldsymbol{F}_{1:R} \bigr). 
    \end{align*}  
If we denote     
  \begin{align*}
\boldsymbol{A} :=  \boldsymbol{F}_{1:R, 1:R}^{-1}   \boldsymbol{F}_{1:R}, \quad  
\boldsymbol{B} := \boldsymbol{F}_{1:R, 1:R}^{-1}   \boldsymbol{F}_{1:R,0},
    \end{align*}
then we can write
\begin{align*}
\boldsymbol{\delta}_0 
   = - \bigl( \boldsymbol{F}_{0,0} -    \boldsymbol{F}_{0,1:R} \boldsymbol{B}
   \bigr)^{-1} 
    \bigl( 
     \boldsymbol{F}_{0} -   \boldsymbol{F}_{0,1:R} \boldsymbol{A}  \bigr), \quad 
 \boldsymbol{\delta}_{1:R} =- \bigl( \boldsymbol{A}+ \boldsymbol{B} \boldsymbol{\delta}_0 \bigr).    
\end{align*}  
Therefore, the system can be solved efficiently as long as
$\boldsymbol{F}_{0,0} -    \boldsymbol{F}_{0,1:R} \boldsymbol{B}$ is invertible.   
Recall that the Newton iteration only requires the computation of $\boldsymbol{f}$ and 
 $\boldsymbol{f}'$, in contrast
with the exact version, which requires the computation of all the derivatives of $f$ up to 
order $R$.

\subsubsection{Computational cost of AIT methods}

From Section~\ref{subsec_newton} we know that for each iteration of Newton's method, we need to compute the vectors     
\begin{align*}
\boldsymbol{\delta}_0 
   = - \bigl( \boldsymbol{F}_{0,0} -    \boldsymbol{F}_{0,1:R} \boldsymbol{B}
   \bigr)^{-1} 
    \bigl( 
     \boldsymbol{F}_{0} -   \boldsymbol{F}_{0,1:R} \boldsymbol{A}  \bigr), \quad 
 \boldsymbol{\delta}_{1:R} =- \bigl( \boldsymbol{A}+ \boldsymbol{B} \boldsymbol{\delta}_0 \bigr),    
\end{align*}  
where $\boldsymbol{A}$ and $\boldsymbol{B}$ are the matrices 
\begin{align*}
\boldsymbol{A} :=  \boldsymbol{F}_{1:R, 1:R}^{-1}   \boldsymbol{F}_{1:R}, \quad  
\boldsymbol{B} := \boldsymbol{F}_{1:R, 1:R}^{-1}   \boldsymbol{F}_{1:R,0}.
    \end{align*}

For the computation of $\boldsymbol{A}$ and $\boldsymbol{B}$, we can exploit  that $\boldsymbol{F}_{1:R,1:R}$ is a blockwise lower triangular matrix with diagonal blocks given by~$-\boldsymbol{I}_M$. Hence, we can obtain for example $\boldsymbol{B}$, by a block forward substitution process, given by 
\begin{equation*}
\boldsymbol{B}_{k}=-\boldsymbol{F}_{k,0}+\sum_{i=1}^{k-1}{\boldsymbol{F}_{k,i}\boldsymbol{B}_{i}},\quad k=1,\dots,R.
\end{equation*}
This algorithm requires $(R^2-R)/{2}$ products of $M\times M$ matrices. 

With respect  to $\boldsymbol{\delta}_0$, the computation of  $\boldsymbol{F}_{0,1:R} \boldsymbol{B}$ requires $R$ products of $M\times M$ matrices, and its $\boldsymbol{LU}$ decomposition requires
$\smash{\bigO(\frac{2}{3}M^3)}$ scalar operations. 
Neglecting operations with lower order cost in $M$, then, we obtain that  the computation of $\boldsymbol{\delta}$ requires
\begin{equation*}C_{\boldsymbol{\delta}} := 
 \left(\frac{R^2+R}{2}+\frac{2}{3}\right)M^3
\end{equation*}
scalar operations. Moreover, 
the formation of the blocks $\boldsymbol{F}_{i,j}$ requires $R^2$~Jacobian matrices of $\boldsymbol{f}$, which yields,  assuming a mean cost per
entry of $\beta$ scalar operations, $R^2M^2\beta$  scalar
operations. In consequence, the (approximate) computational cost of
each Newton iteration for AIT methods is 
\begin{align*} 
 C_{\mathrm{AIT}} =  C_{\boldsymbol{\delta}} +  R^2\beta M^2 = 
 \left(\frac{R^2+R}{2}+\frac{2}{3}\right)M^3+R^2\beta M^2. \end{align*}

\section{Numerical experiments}\label{sec:numexp}

\subsection{Preliminaries} 
In this section, the performance of the AIT methods is analyzed. We first compare 
the AIT methods with their exact counterparts, IT methods, of
the same order.
We have only compared the AIT with the IT methods for scalar
equations since  the implementation for systems of the IT methods is
extremely involved. For linear scalar equations the implementation for
any order is performed using \eqref{eq:pep1234}. 
These methods are compared in terms of error, numerical order {and computational time},
using some scalar problems. 

We then raise two initial-value problems for systems of equations. For those problems, 
the AIT methods are compared with approximate explicit Taylor (AET) methods 
 of the  same order (Baeza et al.\ 2017), 
so as to stress the superior stability of the implicit method. In all the numerical examples we show the numerical errors, computed with $L^1$-norm, and the order of the numerical method, computed by 
$$o(N)=\log_2\left(\left| e(N)/e(N/2)\right|\right),$$
with $e(N)$ standing for the numerical error for $N$ time steps.

\begin{table}[t]\begin{center} 
	\begin{tabular}{|r|ccc|ccc|}	
	\hline
	& \multicolumn{6}{|c|}{$R=2$}\\
	\cline{2-7}
	& \multicolumn{3}{|c|}{IT} &\multicolumn{3}{|c|}{AIT} \\
	\cline{2-7}
	$N$   & $e(N)$   & $o(N)$ &  CPU time  & $e(N)$  & $o(N)$ &  CPU time \\
	\hline
	$10$  & 2.62e-02 &  ---   & 8.5e-02  & 1.38e-02 & ---    & 2.1e-02  \\               
	$20$  & 9.15e-03 & 1.52 & 3.5e-02 & 3.63e-03 & 1.93 & 2.6e-02 \\   
	$40$  & 2.86e-03 & 1.68 & 3.8e-02 & 9.29e-04 & 1.97 & 5.4e-02 \\ 
	$80$  & 8.15e-04 & 1.81 & 4.4e-02 & 2.35e-04 & 1.98 & 1.0e-01 \\  				
	$160$ & 2.19e-04& 1.89 & 5.1e-02 & 5.90e-05 & 1.99 & 2.0e-01 \\               		
	$320$ & 5.70e-05& 1.94 & 4.4e-02 & 1.48e-05 & 2.00 & 4.0e-01 \\              		
	$640$ & 1.45e-05& 1.97 & 3.8e-02 & 3.70e-06 & 2.00 & 7.3e-01 \\              		
	\hline
	\hline
	& \multicolumn{6}{|c|}{$R=3$}\\
	\hline
	$10$   & 1.30e-03 &  ---   & 8.8e-02 & 6.21e-03 &  ---    & 2.9e-02\\               
	$20$   & 2.88e-04 & 2.17 & 4.3e-02 & 9.52e-04 & 2.71 & 3.7e-02\\   
	$40$   & 4.43e-05 & 2.70 & 5.9e-02 & 1.31e-04 & 2.86 & 8.0e-02\\ 
	$80$   & 5.84e-06 & 2.92 & 4.5e-02 & 1.71e-05 & 2.94 & 1.6e-01\\  				
	$160$ & 7.37e-07 & 2.99 & 5.7e-02 & 2.18e-06 & 2.97 & 3.0e-01\\               		
	$320$ & 9.19e-08 & 3.00 & 5.1e-02 & 2.76e-07 & 2.99 & 5.9e-01\\              		
	$640$ & 1.15e-08 & 3.00 & 4.7e-02 & 3.45e-08 & 3.00 & 1.13\\              		
	\hline
	\hline
	& \multicolumn{6}{|c|}{$R=4$}\\
	\hline
	$10$   & 7.55e-04 & ---    & 9.9e-02 & 4.81e-04 & ---     & 5.4e-02\\               
	$20$   & 9.43e-05 & 3.00 & 5.9e-02 & 2.58e-05 & 4.22 & 7.8e-02 \\   
	$40$   & 8.42e-06 & 3.48 & 5.3e-02 & 1.39e-06 & 4.22 & 1.4e-01 \\ 
	$80$   & 6.27e-07 & 3.75 & 6.1e-02 & 7.84e-08 & 4.15 & 2.8e-01  \\  				
	$160$ & 4.26e-08 & 3.88 & 4.6e-02 & 4.61e-09 & 4.09 & 5.4e-01 \\               		
	$320$ & 2.78e-09 & 3.94 & 5.2e-02 & 2.79e-10 & 4.05 &  1.03\\              		
	$640$ & 1.77e-10 & 3.97 & 5.2e-02 & 1.71e-11 & 4.02 &  2.05\\              			
	\hline
	\hline
	& \multicolumn{6}{|c|}{$R=5$}\\
	\hline
	$10$    & 4.32e-05 & ---    & 1.4e-01 & 1.50e-04 & ---    & 8.9e-02\\               
	$20$    & 2.59e-06 & 4.06 & 6.6e-02 & 4.87e-06 & 4.94 & 1.3e-01\\   
	$40$    & 9.73e-08 & 4.73 & 7.3e-02 & 1.54e-07 & 4.98 & 2.2e-01\\ 
	$80$    & 3.14e-09 & 4.95 & 5.5e-02 & 4.86e-09 & 4.99 & 4.3e-01 \\  				
	$160$  & 9.75e-11 & 5.01 & 5.8e-02 & 1.53e-10 & 4.99 & 8.6e-01 \\               		
	$320$  & 3.02e-12 & 5.01 & 5.6e-02 & 4.78e-12 & 5.00 & 1.66\\              		
	$640$  & 9.41e-14 & 5.00 & 6.9e-02 & 1.50e-13 & 5.00 & 3.29\\    
	\hline
	\hline
	& \multicolumn{6}{|c|}{$R=6$}\\
	\cline{1-7}
 	$10$ & 1.59e-05   & ---     &  1.2e-01 & 1.35e-05 & ---      & 1.2e-01\\
 	$20$ & 5.51e-07   &  4.85 &  1.2e-01 & 1.56e-07 & 6.43  & 1.7e-01 \\
 	$40$ & 1.25e-08   &  5.46 &  5.7e-02 & 1.88e-09 & 6.38  & 3.3e-01 \\
 	$80$ & 2.33e-10   &  5.75 &  6.9e-02 & 2.45e-11 & 6.26  & 6.4e-01 \\
 	$160$ & 3.96e-12 &  5.88 &  5.9e-02 & 3.43e-13 & 6.16  &1.30\\
 	$320$ & 6.47e-14 &  5.94 &  7.1e-02 & 5.11e-15 & 6.07  & 2.52\\
 	$640$ & 9.99e-16 &  6.02 &  7.1e-02 & 1.11e-16 & 5.52  & 5.04\\          				
	\cline{1-7}
	\end{tabular}	
	\end{center} 
	\smallskip 	
	\caption{Example~1 (linear scalar problem \eqref{eq:experiment_lineal_escalar}): numerical errors and orders for IT and AIT methods.}
	\label{tab:example_linear}	
\end{table}

	\begin{figure}[t]
	\begin{center} 
				\end{center} 
		\caption{Example~1 (nonlinear scalar problem
                    \eqref{eq:experiment_log_escalar}): performance of the IT and the AIT methods.}
		\label{fig:performance_log}
	\end{figure}

\begin{table}[t]\begin{center} 
	\begin{tabular}{|r|ccc|ccc|}	
	\hline
	& \multicolumn{6}{|c|}{$R=2$} \\
	\cline{2-7}
	& \multicolumn{3}{|c|}{IT} &\multicolumn{3}{|c|}{AIT}\\
	\cline{2-7}
	$N$   & $e(N)$   & $o(N)$ & CPU time & $e(N)$  & $o(N)$  & CPU time  \\
	\hline
	$10$    & 1.21e-03 & ---    & 3.0e-02 & 1.23e-03 & ---     & 4.0e-02 \\               
	$20$    & 2.90e-04 & 2.06 & 6.0e-02 & 2.93e-04 & 2.07 & 4.0e-02 \\   
	$40$    & 7.09e-05 & 2.03 & 6.0e-02 & 7.12e-05 & 2.04 & 4.0e-02 \\ 
	$80$    & 1.75e-05 & 2.02 & 6.0e-02 & 1.76e-05 & 2.02 & 5.0e-02 \\  			
	$160$  & 4.36e-06 & 2.01 & 0.11 & 4.36e-06 & 2.01 & 6.0e-02 \\               		
	$320$  & 1.09e-06 & 2.00 & 0.10 & 1.09e-06 & 2.00 & 0.11 \\              		
	$640$  & 2.71e-07 & 2.00 & 0.16 & 2.71e-07 & 2.00 & 0.21 \\    
	$1280$& 6.77e-08 & 2.00 & 0.26 & 6.78e-08 & 2.00 & 0.37 \\
	$2560$& 1.69e-08 & 2.00 & 0.48 & 1.69e-08 & 2.00 & 0.75 \\          		
	\hline 
	\hline
	& \multicolumn{6}{|c|}{$R=3$} \\
	\hline
	$10$    & 7.52e-05 & ---    & 4.0e-02 & 5.35e-05 & ---     & 4.0e-02  \\               
	$20$    & 8.75e-06 & 3.10 & 6.0e-02 & 5.95e-06 & 3.17 & 4.0e-02\\   
	$40$    & 1.05e-06 & 3.05 & 8.0e-02 & 7.00e-07 & 3.09 & 4.0e-02 \\ 
	$80$    & 1.29e-07 & 3.03 & 9.0e-02 & 8.49e-08 & 3.04 & 5.0e-02\\  			
	$160$  & 1.60e-08 & 3.01 & 0.12 & 1.04e-08 & 3.02 & 9.0e-02\\               		
	$320$  & 1.99e-09 & 3.01 & 0.18 & 1.30e-09 & 3.01 & 0.13\\              		
	$640$  & 2.49e-10 & 3.00 & 0.32 & 1.61e-10 & 3.01 & 0.25\\    
	$1280$& 3.11e-11 & 3.00 & 0.51 & 2.01e-11 & 3.00  & 0.46\\
	$2560$& 3.88e-12 & 3.00 & 0.81 & 2.51e-12 & 3.00 & 0.9\\          		
	\hline
	\hline
	& \multicolumn{6}{|c|}{$R=4$} \\
	\hline
	$10$    & 5.78e-06   &  ---    & 6.0e-02 & 4.93e-06  & ---     & 3.0e-02     \\               
	$20$    & 3.30e-07   & 4.13  & 0.11 & 2.44e-07  & 4.34 & 4.0e-02 \\   
	$40$    & 1.97e-08   & 4.07  & 9.0e-02 & 1.36e-08  & 4.17 & 7.0e-02  \\ 
	$80$    & 1.20e-09   & 4.03  & 0.14 & 8.00e-10  & 4.08 & 7.0e-02   \\  			
	$160$  &  7.43e-11  & 4.02  & 0.23 & 4.86e-11  & 4.04 & 0.19\\               		
	$320$  & 4.62e-12   & 4.01  & 0.76 & 3.00e-12  & 4.02 & 0.19  \\              		
	$640$  & 2.87e-13   & 4.01  & 0.66 & 1.88e-13  & 4.00 & 0.33  \\    
	$1280$& 2.33e-14   & 3.62  & 1.28        & 1.47e-14  & 3.68 & 0.63  \\
	$2560$& 5.66e-15   & 2.04  & 2.17        & 3.77e-15  & 1.96 & 1.27\\ 
	\hline
	\hline   
	&\multicolumn{6}{|c|}{$R=5$} \\
	\hline
	$10$    & 4.52e-07 & ---    & 0.13 & 8.25e-07 & ---      &3.0e-02  \\               
	$20$    & 1.26e-08 & 5.17 & 0.18 & 2.31e-08 & 5.16  &6.0e-02\\   
	$40$    & 3.71e-10 & 5.09 & 0.20 & 6.87e-10 & 5.07  &6.0e-02\\ 
	$80$    & 1.12e-11 & 5.04 & 0.30 & 2.10e-11 & 5.03  &9.0e-02\\  			
	$160$  & 3.42e-13 & 5.04 & 0.48 & 6.53e-13 & 5.01 &0.13\\               		
	$320$  & 4.55e-15 & 6.23 & 0.87 & 2.23e-14 & 4.87 &0.23\\              		
	$640$  & 8.88e-16 & 2.36 & 1.59& 1.55e-15 & 3.84 &0.46\\    
	$1280$& 2.33e-14 & -4.71& 3.11& 2.00e-15 & -0.36 & 0.87\\
	$2560$& 4.22e-15 & 2.47 & 6.09& 2.44e-15 & -0.29 & 1.69\\          	      		
	\hline
	\hline	
	& \multicolumn{6}{|c|}{$R=6$} \\
	\hline
	$10$     & 3.89e-08   &  ---    & 2.04 & 1.52e-07  & ---    & 4.0e-02\\
	$20$     & 5.26e-10   & 6.21  & 1.26 & 1.35e-09  & 6.81 &  6.0e-02  \\
	$40$     & 7.62e-12   & 6.11  & 1.34 & 1.67e-11  & 6.33 &  8.0e-02 \\
	$80$     & 1.17e-13   & 6.03  & 1.50 & 2.19e-13  & 6.26 &  0.11  \\
	$160$   & 4.44e-15   & 4.72  & 1.97 & 1.11e-16  & 10.94 & 0.20  \\
	$320$   & 2.11e-15   & 1.07  & 2.94 & 2.22e-15  & -4.32 & 0.34 \\
	$640$   & 8.88e-16   & 1.25  & 4.92 & 6.66e-16  & 1.74 & 0.64\\
	$1280$ & 2.33e-14   & -4.71 & 8.87& 2.00e-15  & -1.58 &  1.27 \\
	$2560$ & 3.55e-15   & 2.71  & 16.6 & 2.44e-15  & -0.29 & 2.47 \\      				
	\hline
	\end{tabular}
	\end{center} 
	\smallskip 
		\caption{Example~2 (nonlinear scalar problem
                    \eqref{eq:experiment_log_escalar}):  numerical
                    errors, orders {and computational time (in seconds)} for IT and AIT methods.
                }
	\label{tab:example_log}	
\end{table}

	\begin{figure}[t]
	\begin{center} 
		\includegraphics[width=\textwidth,height=0.9\textheight]{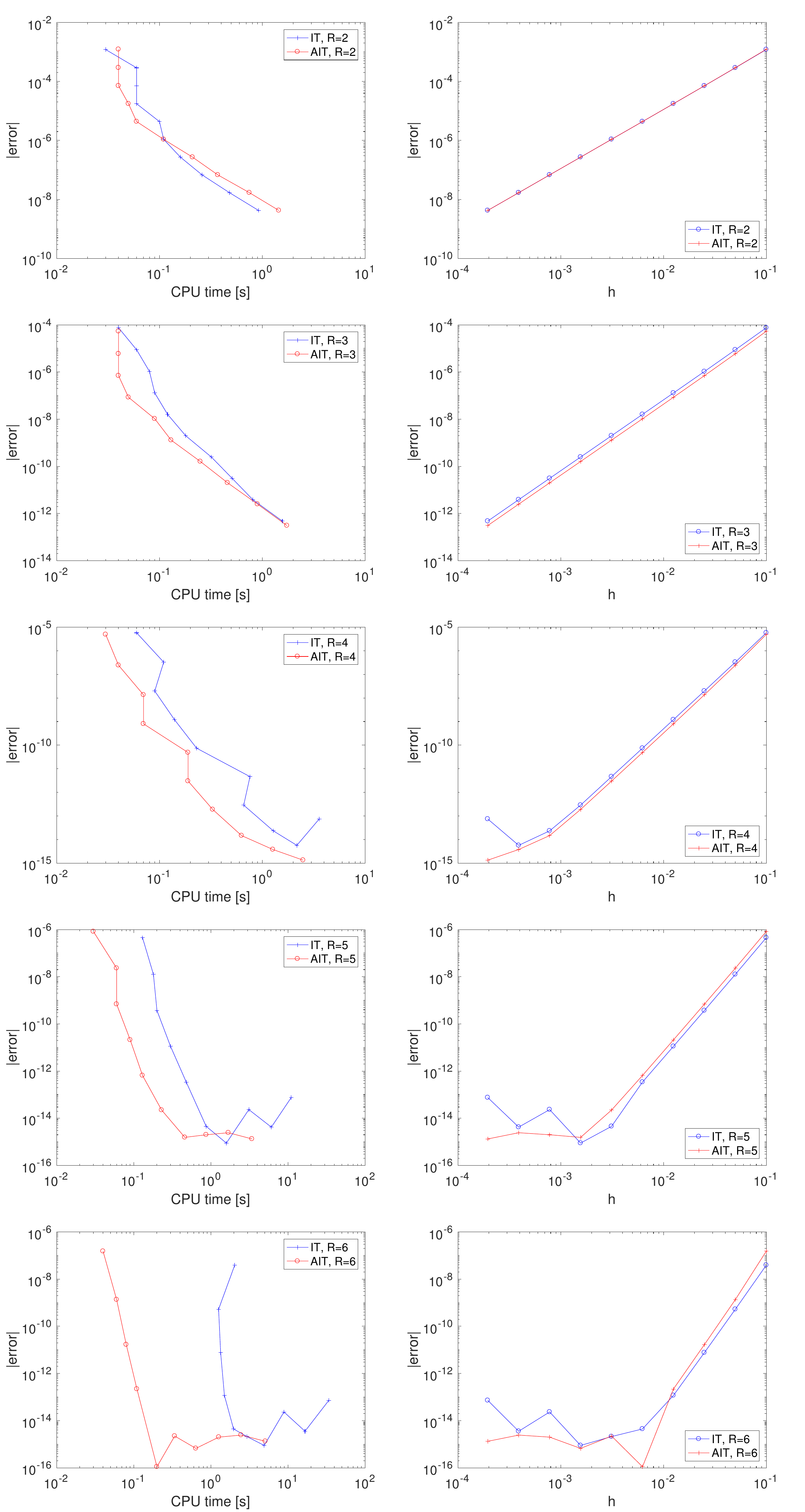}
				\end{center} 
		\caption{Example~2 (nonlinear scalar problem
                    \eqref{eq:experiment_log_escalar}): performance of the IT and the AIT methods.}
		\label{fig:performance_log}
	\end{figure}

\subsection{Examples~1 and~2: scalar equations} 
In Example~1 we consider the linear equation
\begin{equation}\label{eq:experiment_lineal_escalar}
u'= -5u + 5 \sin(2 t) + 2\cos(2 t),\quad u(0)=0,
\end{equation}
with exact solution $u(t) = \sin(2t)$. The results for IT and AIT methods for $T=5$ and orders
$R\in\{2, 3, 4, 5, 6\}$ are collected in Table
\ref{tab:example_linear}, where it can be seen that with both
methods, the expected orders of convergence are recovered in all
cases.  Comparing with the IT methods, we see that the approximate version attains the expected order faster than the exact version, 
but produces a slightly bigger error for coarse resolutions. This fact is possibly due 
to the simplicity of the equation under consideration, which produces a local truncation error smaller than 
the error corresponding to the approximation of derivatives performed in the AIT method and hinders
the correct order of accuracy for the exact method whenever the step size is not small enough.

In Example~2 we  consider the  more involved problem
\begin{equation}\label{eq:experiment_log_escalar}
u'=\log\left(\frac{u+u^3+u^5}{1+u^2+u^4+u^6}\right), \quad u(0)=1,
\end{equation}
and compute its solution up to $T=1$ for orders $R\in\{2,3,4,5,6\}$. The solution computed by 
the AIT method with $R={6}$ and a resolution of  $20000$ points is taken as reference solution.
We can see in Table~\ref{tab:example_log} that the errors for both 
methods are similar and the numerical order converges to the expected values in each case. 
In Figure~\ref{fig:performance_log} we compare the errors obtained by each 
method with respect to the CPU time required to run the algorithm { (left) and with respect the discretization step considered $h$ (right).} It can be seen
that the performance is increasingly favorable to the approximate
method as the order increases, as expected. { Note that, for cases $R= 5, 6$, although the errors are affected by Matlab's computational error, it still can be seen that the AIT method overpowers the IT method in terms of the computational time needed to obtain an approximate solution.}

\subsection{Examples~3 and~4: systems of ODEs}

We consider now two problems modelled by systems of ODEs, used by  Akinfenwa et al.\ (2013)  to test stability properties and accuracy.  Example~3 is a stiff nonlinear problem given by 
\begin{align}\label{eq:experiment_kaps}
\left\{
 \begin{aligned}
y'&= -1002y+1000z^2,\\
z'& = y-z(1+z), \quad t >0; 
\end{aligned} \qquad 
\begin{aligned} 
y(0)=&1,\\
z(0)=&1,
\end{aligned}  
\right.
\end{align}
known as Kaps problem, with exact solution given by
\begin{align*}
y(t)= \mathrm{e}^{-2t},\quad  z(t)= \mathrm{e}^{-t},
\end{align*}
which is independent of the stiffness parameter, $k=-1000$ in this case. We compare the solution
at $T=5$ for the approximate implicit (AIT) and approximate explicit (AET) methods of the same order.
Both schemes recover the expected order, the implicit one achieving it at early stages, see Table~\ref{ex3}. 
Note that  the explicit scheme does not attain good results in terms of accuracy, unless meshes with more than 2000 points are used. The implicit scheme achieves the same error level as the explicit one with meshes with approximately 4 times less nodes.

	\begin{table}[t]\begin{center} 
		\begin{tabular}{|r|cc|cc|cc|cc|cc|cc|}	
			\hline
			& \multicolumn{4}{|c|}{$R=2$} &\multicolumn{4}{|c|}{$R=3$} \\
			\cline{2-9}
		        & \multicolumn{2}{|c|}{AIT} &\multicolumn{2}{|c|}{AET} & \multicolumn{2}{|c|}{AIT} &\multicolumn{2}{|c|}{AET} \\
			\cline{2-9}
			$N$   & $e(N)$   & $o(N)$ & $e(N)$  & $o(N)$  & $e(N)$  & $o(N)$ & $e(N)$ & $o(N)$\\
			\hline
$80$  & 2.12e-05   &1.93 & NaN         &  ---   & 3.31e-07 & 2.92 & NaN         &  ---  \\  				
$160$ & 5.43e-06  &1.96 & NaN         &  ---   & 4.24e-08 & 2.96 & NaN         &  ---  \\               		
$320$ & 1.37e-06  &1.98 &NaN          &  ---   & 5.37e-09 & 2.98 & NaN         &  ---  \\              			
$640$ & 3.45e-07  &1.99 & NaN         &  ---   & 6.76e-10 & 2.99 & NaN         &  ---  \\      
$1280$& 8.66e-08 &1.99 & NaN         &  ---   & 8.47e-11 & 2.99 & NaN          &  ---  \\
$2560$& 2.17e-08 &1.99 & 1.03e-07  & NaN &1.06e-11 & 2.99  & 4.60e-11  & NaN\\ 
$5120$& 5.42e-09 &1.99 & 2.34e-08  & 2.13 &1.32e-12 & 2.99  &  5.75e-12 & 3.00 \\
$10240$&1.35e-09&1.99& 5.84e-09   & 2.00 &1.66e-13 & 2.99  &  7.18e-13 & 3.00   \\    
\hline			
\hline
			& \multicolumn{4}{|c|}{$R=4$} &\multicolumn{4}{|c}{ } \\
\cline{1-5}	
$80$      & 4.13e-09 & 3.92 & NaN        &  ---   & \multicolumn{4}{|c}{ } \\  				
$160$    & 2.65e-10 & 3.96 & NaN        &  ---   & \multicolumn{4}{|c}{ }  \\               		
$320$    & 1.68e-11 & 3.98 & NaN        &  ---   & \multicolumn{4}{|c}{ } \\              		
$640$    & 1.05e-12 & 3.98 & NaN        &  ---   & \multicolumn{4}{|c}{ }  \\      
$1280$  & 6.65e-14 & 3.99 & NaN        &  ---   & \multicolumn{4}{|c}{ }  \\
$2560$  & 4.17e-15 & 3.99 & 1.07e-13 & NaN & \multicolumn{4}{|c}{ }\\ 
$5120$  & 2.70e-16 & 3.94 & 5.65e-15 & 4.24 & \multicolumn{4}{|c}{ } \\
$10240$& 1.95e-17 & 3.79 & 3.33e-16 & 4.09 & \multicolumn{4}{|c}{ }   \\    
\cline{1-5}		
\end{tabular}
\end{center} 
\smallskip 	
\caption{ Example~3 (stiff nonlinear problem \eqref{eq:experiment_kaps}): numerical errors and orders for AIT  and AET methods.}
       \label{ex3}
	\end{table}

	\begin{table}[t]\begin{center} 
		\begin{tabular}{|r|cc|cc|cc|cc|cc|cc|}	
			\hline
			& \multicolumn{4}{|c|}{$R=2$} &\multicolumn{4}{|c|}{$R=3$} \\
			\cline{2-9}
		        & \multicolumn{2}{|c|}{AIT} &\multicolumn{2}{|c|}{AET} & \multicolumn{2}{|c|}{AIT} &\multicolumn{2}{|c|}{AET} \\
			\cline{2-9}
			$N$   & $e(N)$   & $o(N)$ & $e(N)$  & $o(N)$  & $e(N)$  & $o(N)$ & $e(N)$ & $o(N)$\\
			\hline
			$10$  & 5.94e-05 & 2.20  & 3.39e24 &  ---   &       9.59e-06 & 2.45 & 1.24e34   &  ---     \\               
			$20$  & 1.52e-05 & 1.96  & 3.93e37  &  -43.40 &  1.62e-06 & 2.56 & 4.38e50  & -54.97  \\   
			$40$  & 4.10e-06 & 1.89  & 4.11e50  &-43.25    & 2.42e-07 & 2.73 & 4.96e63  &-43.37   \\ 
			$80$  & 1.08e-06 & 1.91  & 2.99e46  &  13.75   & 3.34e-08 & 2.86 & 1.16e46  &58.57 \\  				$160$ & 2.82e-07 & 1.94  & 3.37e-03  &  162.6  & 4.39e-09 & 2.92 & 3.41e-03  &161.2 \\               			$320$ & 7.22e-08 & 1.96  & 7.07e-04  &  2.25   & 5.63e-10 & 2.96 & 2.03e-04  &4.07 \\              			$640$ & 1.82e-08 & 1.98  & 1.67e-04   & 2.08   &  7.12e-11 & 2.98 & 1.95e-05 & 3.38 \\              			\hline
			\hline
			&\multicolumn{4}{|c|}{$R=4$} & \multicolumn{4}{|c|}{$R=5$}\\
			\cline{2-9}
			\hline
			$10$  &  1.69e-06 & 3.05  & 3.29e42   &  ---       & 2.70e-07 & 3.86  & 9.71e49    & --- \\               
			$20$  &  1.56e-07 & 3.43  & 1.43e61   &  -61.91 & 1.28e-08 & 4.39  & 5.80e69& -65.69 \\   
			$40$  &  1.20e-08 & 3.70  & 1.01e72   &-36.05   & 4.97e-10 & 4.69  & 3.26e76 &-22.42 \\ 
			$80$  &  8.32e-10 & 3.85  & 9.17e37   &  113.1   & 1.72e-11 & 4.84  & 4.03e19 &189.0 \\               
			$160$ &  5.48e-11 & 3.92  & 7.99e-04 & 136.4    & 5.69e-13 & 4.92  & 2.66e-04 & 77.01\\               
			$320$ &  3.51e-12 & 3.96  & 3.46e-05 & 4.53      & 1.82e-14 & 4.96  & 5.26e-06  &  5.66 \\
			$640$ &  2.22e-13 & 3.98  & 1.73e-06 & 4.32      & 5.79e-16 & 4.98  & 1.29e-07 & 5.35  \\             
			\hline			
		\end{tabular}
		\end{center} 		\smallskip 		
		\caption{Example~4  (stiff linear problem \eqref{eq:experiment_system_2}): numerical errors and orders for AIT and AET methods.}
		\label{ex4}
	\end{table}

\begin{figure}[t]
	\begin{center} 
		\includegraphics[width=\textwidth,height=0.36\textheight]{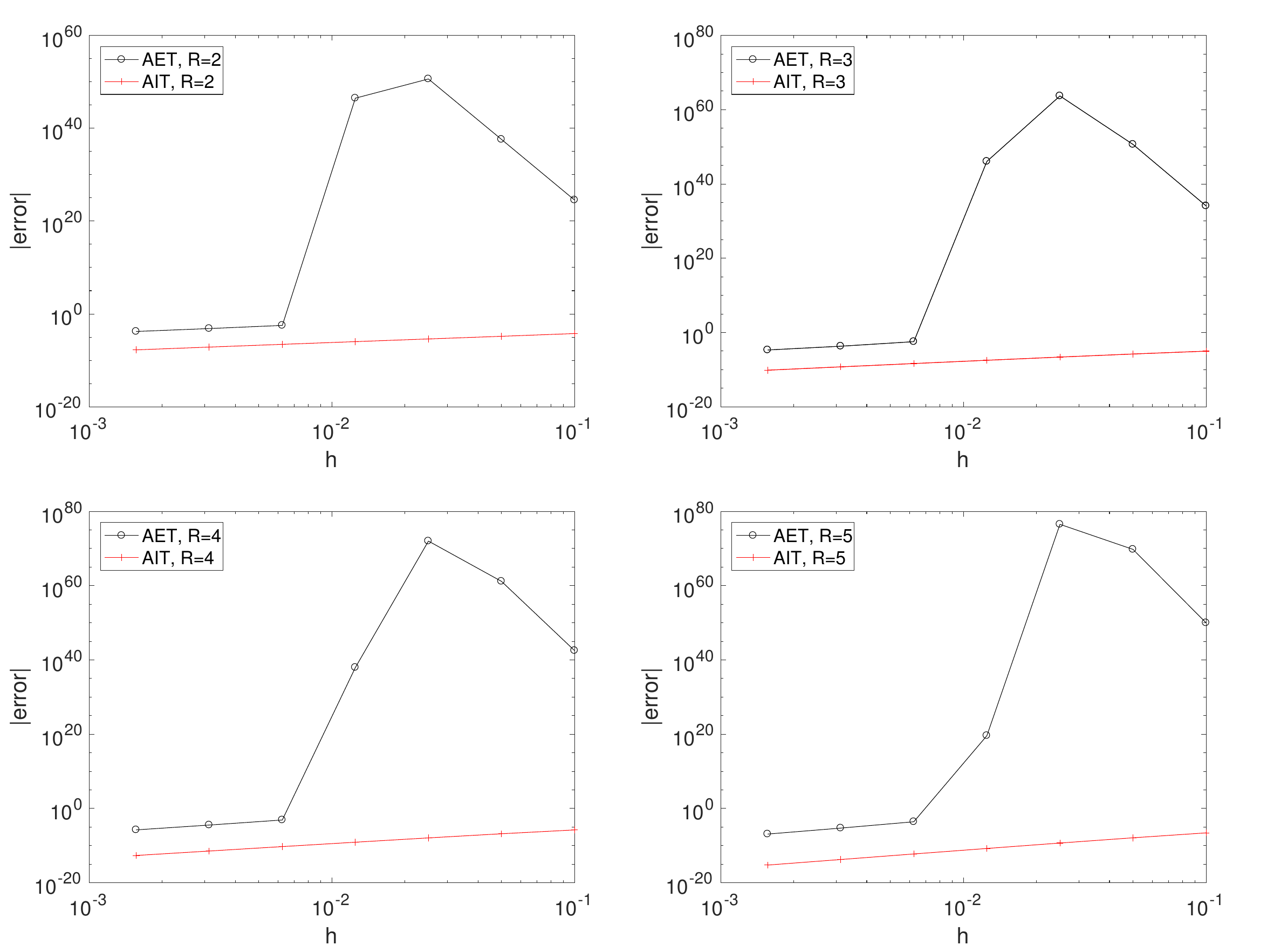}
				\end{center} 
		\caption{Example~4  (stiff linear problem \eqref{eq:experiment_system_2}): performance of the AET and the AIT methods.}
		\label{fig:ex4}
	\end{figure}

Finally, in Example~4 we consider the system of ODEs
\begin{align}\label{eq:experiment_system_2}
\left\{
\begin{aligned}
x'&=-21x+19y-20z,\\
y'&= 19x-21y+20z,\\
z'& =40x-40y-40z,\quad t>0; 
\end{aligned} \qquad 
\begin{aligned} 
x(0)&=1,\\
y(0)&=0,\\
z(0)&=-1,
\end{aligned} 
\right.
\end{align}
also taken from  Akinfenwa et al.\ (2013)  and whose exact solution is given by
\begin{equation*}
\left\{
 \begin{aligned}
x(t)=&\frac12\left(\mathrm{e}^{-2t}+\mathrm{e}^{-40t}(\cos(40t)+\sin(40t))\right),\\
y(t)=&\frac12\left(\mathrm{e}^{-2t}-\mathrm{e}^{-40t}(\cos(40t)+\sin(40t))\right),\\
z(t)=&-\mathrm{e}^{-40t}(\cos(40t)-\sin(40t)).
\end{aligned} 
\right.
\end{equation*}
As in the previous example, the explicit scheme needs more nodes to achieve the same error level as the implicit scheme. For instance, the explicit scheme needs about 64 times more nodes, $N=320$, to obtain the same errors that the implicit scheme attains with $N=5$, see results in Table~\ref{ex4} { and Figure~\ref{fig:ex4}}, proving that the use of the implicit scheme is more appropriate when dealing with stiff problems.

\section{Conclusions}\label{sec:conclusions}

This article is part of ongoing work to develop high-order efficient approximate Taylor ODE solvers. We have reviewed the exact implicit Taylor methods for ODEs and introduced a novel strategy that allows to implement them systematically, although at the cost of differentiating the function in the ODE up to the order of the method.

On the other hand, using the same strategy that led to approximate explicit Taylor methods for ODEs, we define approximate implicit Taylor methods, whose only requirement is the knowledge of function derivatives to build the Jacobian matrix of auxiliary systems of nonlinear equations, to be solved by Newton's method.

While the numerical results essentially confirm that the
  novel approach introduced in this work outperforms the exact version
  in terms of performance, this is not true for low order accuracy, as
  it was expected. The AIT methods are therefore expected to be useful
in the context of ODEs that require to be solved through a
very-high-order implicit scheme.

Finally, it is worth pointing out that it is our purpose
  to further extend this analysis in the context of PDEs, where
  implicit methods are needed in some underlying problems related with
them. High-order methods are being increasingly more demanded to
accurately solve some of these problems, and therefore the AIT methods
may become useful in that context.

\section*{Acknowledgements}
A.B., M.C.M. and P.M. are supported by Spanish MINECO grant  MTM2017-83942-P. P.M. is also supported by  Conicyt/ANID (Chile), project PAI-MEC, folio 80150006 R.B. is supported by  Fondecyt project 1170473;  CRHIAM, Proyecto ANID/Fondap/15130015;  Basal project CONICYT/PIA/AFB170001;  and by the INRIA Associated Team ``Efficient numerical schemes for non-local transport phenomena'' (NOLOCO; 2018--2020),  and   D.Z. is supported by Conicyt/ANID  Fon\-de\-cyt/Post\-doc\-to\-ra\-do/3170077.

\section*{References}

\begin{list}{}{\def\itemindent{-15pt}\def\itemsep{0pt}\def\parsep{0pt}}

\item Abad A, Barrio R, Blesa F, Rodr\'iguez M (2012)  Algorithm 924: TIDES, a Taylor series integrator for differential equation.  ACM Trans Math Softw 39:article~5.

\item  Akinfenwa OA, Jator SN, Yao NM (2013)  Continuous block backward differentiation formula for solving stiff ordinary differential equations. Comput Math Appl 65:996--1005   

\item Baeza A, Boscarino S, Mulet P,  Russo G, Zor\'{\i}o D (2017)  Approximate Taylor methods for ODEs. 
  Comput Fluids 159:156--166

\item Baeza A, Boscarino S, Mulet P,  Russo G, Zor\'{\i}o D  (2020) 
  On the stability of approximate Taylor methods for ODE and the relationship with Runge-Kutta schemes. Preprint, arXiv:1804.03627v1

\item Barrio R,  Rodr\'iguez M, Abad A, Blesa F (2011)  Breaking the limits: The Taylor series method. Appl Math Comput 217:7940--7954 

\item Dennis Jr.\ JE,   Schnabel RB (1996)  
   Numerical Methods for Unconstrained Optimization and Nonlinear
 Equations.   Classics in Applied Mathematics vol.~16,  SIAM,  Philadelphia

\item  Fa\`a di Bruno F (1855) Sullo sviluppo delle funzioni. Annali di Scienze Matematiche e Fisiche 
6:479--480   

\item Hairer E, Wanner G (1996) ,Solving Ordinary Differential Equations II. Stiff and Differential-Algebraic Problems. 2nd edition. Springer Series in Comput Math vol. 14.

\item  Jorba \`A,  Zou M (2005) A software package for the numerical integration of ODEs by means of high-order Taylor methods. Exp Math 14:99--117

\item Kirlinger G, Corliss GF (1991) On implicit Taylor series methods for stiff ODEs, Argonne National Laboratory technical report ANL/CP-74795 

\item Miletics E, Moln\'arka G (2004) 
 Taylor series method with numerical derivatives for numerical solution of ODE initial value problems.
   J Comput Methods Sci Eng. 4:105--114
   
\item Miletics E, Moln\'arka G (2005)
 Implicit extension of Taylor series method with numerical derivatives for initial value problems. 
   Comput Math Appl 50:1167--1177
   
 \item Qiu J,  Shu CW (2003)  Finite difference WENO schemes with Lax-Wendroff-type time discretizations. SIAM J Sci Comput  24:2185--2198 

\item Scott JR (2000)   Solving ODE initial value problems with implicit Taylor series methods, NASA technical memorandum TM-2000-209400

\item Zor\'{\i}o D,  Baeza A,  Mulet P (2017)  An approximate Lax-Wendroff-type procedure for high order accurate schemes for hyperbolic conservation laws. J Sci Comput 71:246--273

\end{list}

\end{document}